\documentclass[12pt,reqno]{amsart}

\usepackage{amssymb}
\usepackage{amsthm}
\usepackage{amsmath}
\usepackage{amsxtra}
\usepackage{mathrsfs}
\usepackage{enumerate}
\usepackage{subfig}
\usepackage{tikz}
\usepackage{tikz-cd}
\newtheorem{theorem}{Theorem}[section]

\theoremstyle{remark}

\numberwithin{table}{section}

\newcommand{\C}{\mathbb C}

\newcommand{\Q}{\mathbb Q}

\newcommand{\Z}{\mathbb Z}

\newcommand{\slp}{SL_2(\Z/p \Z)}
\newcommand{\slpsq}{SL_2(\Z/p^2 \Z)}
\newcommand{\slpn}{\text{SL}_2(\Z/p^n \Z)}

\renewcommand{\arraystretch}{1.5}

\usepackage{ wasysym }

\usepackage{fullpage}
\usepackage[all,cmtip]{xy}

\begin{document}

\title[The Computation of Fourier Transforms on $SL_2(\Z/p^n \Z)$]{The Computation of Fourier Transforms on $SL_2(\Z/p^n \Z)$ and Related Numerical Experiments}

\author{Benjamin K. Breen}
\address{Department of Mathematics, Dartmouth College, Hanover, NH 03755}
\email{Benjamin.K.Breen.GR@dartmouth.edu}

\author{Daryl R. Deford}
\address{Department of Mathematics, Dartmouth College, Hanover, NH 03755}
\email{ddeford@math.dartmouth.edu}

\author{Jason D. Linehan}
\address{Department of Mathematics, Dartmouth College, Hanover, NH 03755}
\email{Jason.D.Linehan.GR@dartmouth.edu}

\author{Daniel N. Rockmore}
\address{Departments of Mathematics and Computer Science, Dartmouth College, Hanover, NH 03755}
\email{rockmore@math.dartmouth.edu}

\subjclass[2000]{To be filled in.}

\date{\today}



\begin{abstract}
We detail an explicit construction 
of ordinary irreducible representations for the family of finite groups $SL_2(\Z /p^n \Z)$ for odd primes $p$ and $n\geq 2$. For $n=2$, the construction is a complete set of irreducible complex representations, while for $n>2$, all but a handful are obtained. We also produce  an algorithm for the computation of a Fourier transform for a function on $SL_2(\Z /p^2 \Z)$.  With this in hand we explore the spectrum of a collection of Cayley graphs on these groups, extending analogous computations for Cayley graphs on $SL_2(\Z/p \Z)$ and suggesting conjectures for the expansion properties of such graphs. 
\end{abstract}

\maketitle

\null\hfill\begin{tabular}[t]{l@{}}
\text{} \\
\end{tabular}

\section{Introduction}

Ever since the foundational work of Lubtozky, Phillips, and Sarnak \cite{LPS} the expansion properties of the groups $SL_2(\Z /p \Z)$ have been of significant interest. Work of Lafferty and Rockmore \cite{LR92,LR93}  conjectured that random sets of generators for these groups (suitably defined) would generically give rise to expander graphs, a conjecture ultimately settled by Bourgain and Gamburd \cite{BG}. Lafferty and Rockmore studied these Cayley graphs by computing the Fourier transform of the characteristic function of the generators at a complete set of irreducible representations for these groups \cite{LR99}. Since the dimension of the individual Fourier transforms is much smaller than the order of the group (in this case, roughly $p$ versus $p^3$) more extensive calculations could be accomplished. Additionally, this ``microanalysis" suggested more fine scale conjectures about the spectra, that  not only  the expansion properties of random generators were generic, but even more, that the spectra of individual Fourier transforms was generic. This conjecture was only recently settled by Rivin and Sardari \cite{RS}. 

Analogous questions for the related family of groups $SL_2(\Z/p^n \Z)$ are natural generalizations and have yet to be investigated. The goal of performing the necessary numerical experiments, \`{a} la \cite{LR92,LR93,LR99} are the motivation for this paper. These experiments can only be accomplished if we have at hand a way to efficiently and explicitly construct a complete set of irreducible complex matrix representations for these groups. Achieving this first goal has the added attraction of producing yet another infinite family of important finite groups whose irreducible representations are then explicitly accessible. This in turn produces the basic material for considering the problem of efficiently computing a Fourier transform on these groups, thereby relating this work to another broad field of current research (see e.g., \cite{MaslenRockmoreNotices,sovII}). 

Our construction follows and clarifies early work of Tanaka \cite{Tanaka67-2,Tanaka67-1} wherein the irreducible representations were first written down at varying levels of detail. This is the contents of Section 2. In Section 3 we give an example illustrating the construction of Section 2. In Section 4 we show a number of numerical experiments on Cayley graphs of $SL_2(\Z/p^n \Z)$. We conclude in Section 5 by outlining possible future work. Sage code implementing this methods is available at \cite{GH_Code}.

\section{Construction of the irreducible representations of $SL_2(\Z/p^n \Z)$ for $n \geq 2$}

In this section we give an explicit construction of the irreducible representations for $SL_2(\Z/p^n \Z)$ for $n \geq 2$. Our work derives from that of Tanaka  \cite{Tanaka67-1},  filling in and correcting details as necessary as we work toward the goal of computing Fourier transforms on these groups. 

The representations sort immediately into those representations that either don't or do come from representations of $SL_2(\Z/p^{n-1} \Z)$ via the composition

\begin{equation}\label{eq:proj}
SL_2(\Z/p^n \Z) \;\; {\overset{\pi }{ \longrightarrow}} \;\; SL_2(\Z/p^{n-1} \Z) {\overset{\rho }{ \longrightarrow}}  \;\; GL_d(\C)
\end{equation}
where $\pi$ is the natural projection map (via reduction of the matrix entries mod $p^{n-1}$) from $SL_2(\Z/p^n \Z)$ to $ SL_2(\Z/p^{n-1} \Z) $ and $\rho$ is any representation of $SL_2(\Z /p^{n-1} \Z)$. We will call these {\em (irreducible) quotient representations.} This gives a recursive roadmap for the construction: To compute the irreducible representations of  $SL_2(\Z/p^n \Z)$, first compute all the irreducible representations of $SL_2(\Z/p^{n-1} \Z)$ and then compute any representations that don't arise in this way. Continuing down this chain of projections we arise at a ``base case" of $n=1$ and the consideration of the irreducible representations of $SL_2(\Z/p \Z)$. For that we draw on the explicit construction in \cite{LR92}
.

Tanaka  \cite{Tanaka67-1} outlines the construction of the non-quotient representations.  
Some readers may note that the mathematical structure underlying this construction is derived from considering quadratic extensions of $p$-adic fields and some of the language used reflects this fact.
 

\subsection{Outline of Construction}

Following  \cite{Tanaka67-1} we produce several large (and reducible) representations of $SL_2(\Z/p^n \Z)$ that depend on three parameters: $k$, $\Delta$, $\sigma$. By varying these parameters,  we produce representations that contain all non-quotient irreducible representations of $SL_2(\Z/p^n \Z)$.  We then extract the individual irreducible representations by computing the actions on carefully chosen subspaces of the reducible space. For the case of $n=2$ we obtain explicit direct constructions of all the irreducible representations. For $n>2$, we obtain (in some strong sense made clear below) almost all the individual irreducibles, leaving a few as contained within some larger reducible representation. We hope in future work to be able to isolate the few remaining holdouts. 

We first take up the parameter $k$, which ranges $0\leq k \leq n$. With $k$ fixed we consider the group $G = \Z/p^n \Z \times \Z/p^{n -k}  \Z$  and the associated free complex vector space,  $\C [G]$, with canonical basis given by the group elements.\footnote{This is also the notation used for the group algebra generated by  the additive group $G$, but we will not need this additional structure.}   When $k=n$ we have the familiar additive group $\Z/p^n$. For $k<n$ more work is needed. An action of $SL_2(\Z/p^n \Z)$ of $\C[G]$ is then defined that requires a choice of the parameter $\sigma$ as well as a  ring structure on $G$ that depends on the parameter $\Delta$. (For ease of reference, we keep most of Tanaka's original notation). This ring we denote as $G(k,\Delta)$. The resulting representation given by turning this ring into a complex vector space indexed by the group (ring) elements is denoted $R_k(\Delta, \sigma)$.

The representation $R_k(\Delta, \sigma)$ is reducible and the most difficult part of the construction comes from finding the irreducible representations inside $R_k(\Delta, \sigma)$. This is accomplished by finding a specified (multiplicative) abelian group $C \subset G(k,\Delta)$. A specific set of so-called {\em principal} characters $\chi \in \widehat{C}$ then determine a basis for an irreducible subrepresentation inside $R_k(\Delta, \sigma)$. We denote  these irreducible sub-representations by $R_k(\Delta, \sigma, \chi)$. 

To explicitly specify the representation, we only need to work out the action of generators  of $SL_2(\Z/p^n \Z)$. These are given by the factors in the {\em Bruhat decomposition}.  For this, let $D$ denote the diagonal subgroup of $SL_2(\Z/p^n \Z)$ and $U$, the {\em unipotent subgroup}. If we use the notation

\[ d_a = \left(
\begin{array}{cc}
a & 0 \\
0 & a^{-1} \\
\end{array} \right) \quad  \quad  u_b = \left(
\begin{array}{cc}
1 & b \\
0 & 1 \\
\end{array} \right)  \quad \quad w = \left(
\begin{array}{cc}
0 & -1 \\
1 & 0 \\
\end{array} \right)\]
where $a \in (\Z/p^n \Z)^\times$ and $b \in \Z/p^n \Z$, then given   

\[ A =  \left( \begin{array}{cc}
\alpha & \beta \\
\gamma & \delta \\
\end{array} \right)\in \text{SL$_2$}(\Z/p^n \Z)\]
its {\em Bruhat decomposition} is given by
\begin{equation} A =\begin{cases} u_{\alpha\gamma^{-1}}\cdot w \cdot d_\gamma \cdot u_{\delta \gamma^{-1}}  & \text{if} \quad \gamma \not \equiv 0  \pmod{p}\\
w \cdot u_{-\gamma\alpha^{-1}} \cdot w \cdot  d_{-\alpha}\cdot  u_{\beta\alpha^{-1}}  & \text{if} \left\{\begin{array}{l}  \gamma \equiv 0 \pmod{p}\\  \alpha \not \equiv 0  \pmod{p}\end{array}\right. . \end{cases}\end{equation}\label{Bruhat}

From the Bruhat decomposition it follows immediately that $|SL_2(\Z/p^n \Z)| = p^np^n(p^n-p^{n-1}) + p^np^{n-1}(p^n-p^{n-1})$.  And to compute the matrix representation of any group element it suffices to specify (compute) the actions of $d_a$, $u_b$, and $w$ on a given representation space. The details now follow.

\subsection{The Ring $G(k,\Delta)$: The case of $k<n$.}


For $0\leq k\leq n-1$, let  $G$ denote the abelian (additive) group (the dependence on $k$ is understood) with point set

\[ G \: \:  = \: \:  \Z/p^n \Z \times \Z/p^{n -k}  \Z. \]

As usual, let $\C[G]$ denote the free $\C$ vector space  with basis $G$  $$ \{\alpha = \sum_{g \in G} \alpha_g g:\alpha_g\in\C\}$$ which will be the representation space for each of our representations.

Let $\Delta'$  be a square-free positive integer, such that $(\Delta',p)=1$ and define $\Delta=p^k\Delta'$. We  put a ring structure on $G$, by defining a multiplication determined by the parameter $\Delta'$ and denote the subsequent object $G(k,\Delta)$. To define the multiplication structure we make use of the embedding $\Z/p^{n-k}\Z\hookrightarrow\Z/p^n\Z$ given by $a \mapsto p^k \cdot a$. This gives a map from elements of $\Z/p^{n-k}\Z$ into $\Z/p^n\Z$, for example if $g_2 \in \Z/p^{n-k} \Z$ then $g_2\Delta = p^k \cdot g_2 \Delta' \in \Z/p^n\Z$. For any $g = (g_1,g_2) \in G$, let $\widetilde{g} = g_1 + g_2 \sqrt{-\Delta}$ and define a multiplication $*$ on $G$ by 

\begin{align*}
g*h :=&\quad \widetilde{u} \cdot \widetilde{v} \\
 =&\quad \left(g_1 + g_2 \sqrt{-\Delta}\right) \: \:  \cdot \: \: \left(h_1 +  h_2 \sqrt{-\Delta} \right)\\
  \equiv&\quad y_1 +  y_2\sqrt{-\Delta} \\ 
  =& \quad \widetilde{y} \\ 
  \end{align*}

\noindent where $y_1 = g_1h_1 - (g_2h_2) \Delta \pmod{p^n}$ and $ y_2 = g_1h_2 + g_2h_1 \pmod{p^{n-k}}$.



Having endowed $(G(\Delta,k),*)$ with a multiplicative structure, we further define  {\em conjugation} (${\bar{\hspace{1em}} }$), {\em norm} (Nm) and {\em trace} (Tr) maps:
$$\begin{array}{lll}
{\bar{\hspace{1em}} }: G   \to  G &  \qquad \qquad \overline{(g_1,g_2)} \: \:  \mapsto \: \:   (g_1,-g_2) \\
\vspace{-.2cm}
 \mbox{Nm:}  \; G   \to    \Z/p^n \Z &  \qquad \qquad  \text{Nm}(g_1,g_2) \: \:  \mapsto \: \:   g_1^2 + \Delta g_2^2  & \text{(Multiplicative)} \\
\vspace{-.2cm}
 \text{Tr:} \; G   \to  \Z/p^n \Z & \qquad \qquad  \text{Tr}(g_1,g_2) \: \:  \mapsto \: \:   2g_1 \: \: \:  & \text{(Additive)} \end{array}$$


\medskip

Note that Nm and Tr are multiplicative and additive with respect to group multiplication on $G(k,\Delta)$ and the standard componentwise additive structure on $G$: $\text{Tr}(g + h) = \text{Tr}(g ) + 
\text{Tr}( h)$ and $\text{Nm}(g*h) = \text{Nm}(g)\text{Nm}(h)$. We use the same embedding of $\Z/p^{n-k}\Z\hookrightarrow\Z/p^n\Z$ to compute the norm. 



\subsection{The representation  $R_k(\Delta, \sigma)$: The case of $k<n$.}\label{Representation}

We define a representation space $R_k(\Delta, \sigma)$  with map $T_\sigma$ by specifying the actions of the Bruhat elements $d_a, u_b$ and $w$ on $\C[G]$. We first let $\sigma \in \Z$ determine a $p^n$th-root of unity $\zeta_{\sigma} = e^{2 \pi i \sigma / p^n}$ for the characters. Now define the representation $R_k(\Delta, \sigma)$ via the actions 


\begin{equation}\label{eq:a}T_\sigma(d_a) \left [ \: g \:  \right ] = \left ( \frac{a}{p} \right )^k {g \cdot a^{-1}}\end{equation}


\begin{equation}\label{eq:b} T_\sigma(u_b) \left [ \: g \:  \right ] = \zeta_{\sigma}^{\: b Nm(g)} g \end{equation}


\begin{equation}\label{eq:w}T_\sigma(w) \left [ \: g \:  \right ] =  c \displaystyle \sum_{h \in G} \zeta_{\sigma}^{\: -Tr(g *\overline{h})}  h  \end{equation}
where  $\left ( \frac{a}{p} \right )$ denotes the Legendre symbol, and if $g = (g_1,g_2)$ then $g \cdot a =(ag_1,ag_2)$ and $c$ is a constant 
given by

\[ c = p^{-n +(k/2)} \left ( \frac{\Delta'}{p} \right )^{n -k} \left ( \frac{\sigma}{p} \right )^{k} \cdot e \]
where 

\[e = 
\begin{cases} 1 & k \text{ odd} \hspace{1cm} \left ( \frac{-1}{p} \right ) = 1 \\
-i & k \text{ odd} \hspace{1cm} \left ( \frac{-1}{p} \right ) = -1 \\
-1^n & k  \text{ even} \end{cases}  \]

Tanaka \cite{Tanaka67-2} defines these actions on the dual space of functions. In order  to effectively implement these methods, we select the standard basis to make use of the non-canonical isomorphism of $\C[G]^*$ and $\C[G]$.

\subsection{Irreducible sub-representations $R_k(\chi, \Delta, \sigma)$: The case $k<n$.}

Having built the large representation 
$R_k(\Delta, \sigma)$ our next goal is to carve out irreducible subspaces of this large representation. To do this we consider the monoid $(G, *)$ under multiplication 
 and define

\begin{equation}
C \: = \: \{ \: u \in G \: \: : \: \:  \text{Nm}(u) \equiv 1 \pmod{p^n} \:  \} 
\end{equation}\label{eq:Cdef}
which is an abelian group within the monoid. Additionally, we define
\[ C_{\: n-1}  =  \{ \: c \in C  :  \quad u_1 \equiv 1  \pmod{p^{n-1}} \quad \text{and} \quad u_2 \equiv 0   \pmod{p^{n-1-k}} \: \}. \]


We say a character $\chi$ on $C$, is \textbf{principal} if the restriction of $\chi$ to $C_{n-1}$ is nontrivial. Otherwise we say it is {\bf decomposable.} All the non-quotient irreducibles can be realized in terms of principal characters. For a principal character $\chi$, we consider the induced vector space $\operatorname{Ind}(V_\chi)$:
\[ \operatorname{Ind}(V_\chi) :=  \{   \alpha \in \C[G]  \: \: : \: \:  \alpha_{c\: * \: g} = \chi(c) \alpha_{g} \quad \quad c \in C , g \in G \: \}. \]

When $k<n$, the restriction of the representation on $R_k(\Delta, \sigma)$ given in Section \ref{Representation} to the subspace $\operatorname{Ind}(V_\chi)$ is irreducible (cf. \cite{Tanaka67-2} Section 4, Theorem i)).
 We denote this irreducible representation as $R_k(\chi, \Delta, \sigma)$. 
 
 \subsubsection{Final representations: $R_k(\Delta,\sigma)$ for $k=n$} \label{k=n} 

The above construction produces all the irreducibles for the cases $k<n$ with Theorem~\ref{thm:sl2smallk} summarizing the contributions in the case of $n=2$. In this section we show how the final set of non-quotient irreducibles, found in $R_n(\Delta,\sigma)$, are determined. This is the case of $k=n$ and the above construction has to be altered slightly in order to find irreducible representations and depends on   \cite{Kloosterman46} . Following the construction as in the case of $k<n$ we see that

\[ G \cong \Z/p^n\Z \times \Z/p^{n-n}\Z  \cong \Z/p^n\Z \times 0 \cong \Z/p^n\Z. \]

Thus,  the imposed multiplication structure no longer depends on $\Delta$ and is just the usual multiplication on $\Z/p^n\Z$. Note that this also implies that the representation only depends on $\sigma$, so we will denote it as $R_n(\sigma)$. 

Within $R_n(\sigma)$ there are two important subrepresentations $R_n(\chi_i,\sigma)$ for $i =\pm 1$ that we need to identify. To find them, observe first that in this case, the subgroup $C$ reduces to the two element group: 
\[ C = \{ (1,0), \: (p^n-1,0) \}. \]
The corresponding character group $\widehat{C}$ has two elements: $\chi_1$ the trivial character and $\chi_{-1}$ the non-trivial character. Then according to \cite{Kloosterman46} (pp. 371-372) we have the following:

\begin{theorem}
The representation  $R_n(\sigma)$  has a decomposition into 
$$R_n(\sigma) = R_n(\chi_1, \sigma) \oplus R_n(\chi_{-1}, \sigma)$$
with 
$${\text{dim}}_{\C} R_n(\chi_1, \sigma) = (p^n+1)/2  {\text{ and  }} {\text{dim}}_{\C} R_n(\chi_{-1}, \sigma) = (p^n-1)/2.$$
\end{theorem}


In contrast to the previous situation, the representations $R_n(\chi_1, \sigma)$ and $R_n(\chi_{-1}, \sigma)$  are reducible when $n \geq 2$ and $n \geq 3$ respectively.    To see this, let $\pi$ denote the natural projection $\Z/p^n\Z  \;\; {\overset{\pi }{ \longrightarrow}} \;\;  \Z/p^{n-1}\Z $ and $i$ the inclusion $\Z/p^n\Z  \;\; {\overset{i }{ \longrightarrow}} \;\;  p\Z/p^{n+1}\Z $. The composition of $i$ and the inverse image $\pi^{-1}$ yields a map 

\[ \xymatrix@1{ \Z/p^{n}\Z  \ar[dr]_{\pi} & & \ar@{-->}[ll]_{\pi^{-1}\circ i} \ar[dl]^{i} \Z/p^{n-2}\Z \\ & \Z/p^{n-1}\Z & }\]

\noindent expressed as $\pi^{-1}\circ i : \Z/p^{n-2}\Z \to  \Z/p^{n}\Z$ which induces a map on the group algebras.

\[ \C[\Z/p^{n-2}\Z]  \;\; {\overset{\pi^{-1}\circ i}{ \longrightarrow}} \;\; \C[\Z/p^{n}\Z] \]



This is actually an embedding of the representation $R_{n-2}(\sigma) \hookrightarrow R_{n}(\sigma)$ and it actually carries $R_{n-2}(\chi_1,\sigma) \hookrightarrow R_{n}(\chi_1, \sigma)$ and $R_{n-2}(\chi_{-1},\sigma) \hookrightarrow R_{n}(\chi_{-1}, \sigma)$. Thus the representation $R_n(\sigma)$ decomposes as 

\[ R_n(\sigma) \: \;  \cong \: \: R_{n}(\chi_1, \sigma) \oplus R_{n}(\chi_{-1}, \sigma)\]

\[  \cong  \: \: \left [ R_{n}(\chi_1, \sigma) - R_{n-2}(\chi_1, \sigma) \right ] \: \: \oplus \: \:  \left [ R_{n}(\chi_{-1}, \sigma) - R_{n-2}(\chi_{-1}, \sigma) \right ] \: \: \oplus \: \: R_{n-2}(\sigma) \]

\[  \cong  \: \: R_{n}(\chi_1, \sigma)_1 \: \: \oplus \: \:   R_{n}(\chi_{-1}, \sigma)_1 \: \: \oplus \: \: R_{n-2}(\sigma) \]

The last set of irreducible, non-quotient representations is actually $R_{n}(\chi_i, \sigma)_1$ for $i = \pm1$. A basis [\cite{nobs1976irreduziblen} pg. 509 Satz 4] is given by $B_0(\pm) \cup B_1(\pm)$ where

\begin{enumerate}
\item  $ B_0(\pm) = \left\{ \quad (x,0)\pm(-x,0) \quad : \quad x \in (\Z/p^n\Z)^\times \quad  \textrm{and} \quad  1\leq x\leq \dfrac{p^n-1}2\right\}$ 
\item $B_1(\pm) = \left\{ 
\displaystyle{\sum_{a \: \in \: [0 \dots p]}} \zeta^{ka} [ (py +ap^{n-1}) \pm (-py -ap^{n-1}) ] \:  :  \begin{array}{l} 0 \leq k \leq \dfrac{p-1}2   \: \: \textrm{and}\: \ \\ 
 0 \leq y \leq p^{n-2}-1 \end{array}  \right\}.$
\end{enumerate}

where $\zeta = e^{2\pi i/ p}$.

 \subsection{Summary of construction}
 Although the construction described above is fairly complex,  we can summarize the discussion as a description of the representations according to the values of the three parameters $k,\Delta ,$ and $\sigma$:

$$\begin{array}{lll}
k: & \mbox{Determines the space for the representation.} &  ( 0 \leq k \leq n).\\
\Delta: &  \mbox{Determines arithmetic, depends on }k.  &  \left\{\begin{array}{l} \Delta = \Delta' p^k \in \Z \mbox{ with }\\  (\Delta', p) =1 \mbox{ and squarefree.}\end{array}\right. \\
\chi: & \mbox{Non-decomposable character on an abelian group } C. & \mbox{ Depends on }\Delta .\\
\sigma: &  \mbox{Used to construct characters. } & \sigma \in \Z ;   \qquad (\sigma, p) = 1.
\end{array}$$


The results of Tanaka \cite{Tanaka67-1} specify which choices of parameters lead to inequivalent representations and describe the number and dimension of each space. We record these results in the following theorem.

\begin{theorem}[Tanaka]
The equivalence of representations $R_k(\chi,\Delta,\sigma)$ depends on the value $k$ as follows:

\begin{itemize}
\item[(i)] When $k=0$ two representations $R_0(\chi_1,\Delta_1,\sigma_1)$ and  $R_0(\chi_2,\Delta_2,\sigma_2)$ are equivalent if and only if $\left ( \frac{\Delta_1}{p} \right ) = \left ( \frac{\Delta_2}{p} \right )$ and $\chi_1=\chi_2$ or $\chi_1=\chi_2^{-1}$.

\item[(iii)]  When $1\leq k\leq n-1$  two representations $R_k(\chi_1,\Delta_1,\sigma_1)$ and  $R_k(\chi_2,\Delta_2,\sigma_2)$ are equivalent if and only if $\left ( \frac{\Delta_1}{p} \right ) = \left ( \frac{\Delta_2}{p} \right )$, $\left ( \frac{\sigma_1}{p} \right ) = \left ( \frac{\sigma_2}{p} \right )$ and $\chi_1=\chi_2$ or $\chi_1=\chi_2^{-1}$.

\item[(iii)]   When $k=n$ two representations $R_n(\chi_1,\Delta_1,\sigma_1)$ and  $R_n(\chi_2,\Delta_2,\sigma_2)$ are equivalent if and only if $\left ( \frac{\sigma_1}{p} \right ) = \left ( \frac{\sigma_2}{p} \right )$ and $\chi_1 = \chi_2$.
\end{itemize}
Table~2.1 below summarizes the properties of the objects constructed above:
\begin{table}[!h]
\begin{tabular}{|l|c|c|c|}
\hline
$R_k(\chi,\Delta,\sigma)$&$k=0$&$1\leq k\leq n-1$&$k=n$\\
\hline
Number&$(p^{n-1})(p-1)$&$4(p^{n-k}-p^{n-k-1})$&$4$\\
\hline
Dimension&$p^n+  ( \frac{-\Delta'}{p} )p^{n-1}$&$p^{n-2}(p^2-1)/2$&$p^{n-2}(p^2-1)/2$\\
\hline
$|C|$&$p^{n-1}(p- (\frac{-\Delta'}{p} ))$&$2p^{n-k}$&$2$\\ 
\hline
$|\chi|$&$p^{n-2}(p-1) (p- (\frac{-\Delta'}{p} ) )$&$2(p^{n-k}-p^{n-k-1})$&$2$\\
\hline

\end{tabular}\label{tbl:pnreps}
\caption{Classification of non-quotient irreducible representations of $SL_2(\Z/p^n\Z)$.}

\end{table}
\end{theorem}

\subsection{The irreducible representations of $SL_2(\Z/p^2\Z)$.}

For $SL_2(\Z/p^2\Z)$ we can complete the story explicitly. In this case the group $C$ (cf. Equation~\ref{eq:Cdef}) is always cyclic \cite{Kloosterman46}. This has numerous computational advantages, both for determining the principal characters and evaluating equations (3), (4), and (5). In fact, choices of the characters in this case can be computed directly for all parameter choices. We record these characters in the following theorem, noting that $\chi_j$ and $\chi_\ell$ give equivalent representations if and only if $\chi_j=\chi_\ell$ or $\chi_j=-\chi_\ell$ \cite{Tanaka67-1}. 
 
 \begin{theorem}\label{thm:sl2smallk}
In the case $n=2$ a complete set of principal characters 
is given by
 \begin{itemize}
 \item When $k=0$ and $\left(\dfrac{-\Delta'}{p}\right)=1$   $$\{ 0< \ell <p(p-1)\ :\  (\ell,p)=1\}.$$
 \item When $k=0$ and $\left(\dfrac{-\Delta'}{p}\right)=-1$  $$\{ 0< \ell <p(p+1)\ :\ (\ell,p)=1\}.$$
 \item When $k=1$  $$\{0<\ell<2p\ :\ (\ell,p)=1\}.$$
 \end{itemize}
 \end{theorem}


All that is left is to consider $k=n = 2$. But in this  case the representation $R_2(\chi_{-1}, \sigma)$ is irreducible of dimension $(p^2-1)/2$. However, the representation $R_2(\chi_1, \sigma)$ is reducible and contains a copy of the trivial representation with the remaining representation $R'_2(\chi_1, \sigma)$ of dimension $(p^2-1)/2$ irreducible. Among our new results is an explicit basis for computing these representations: 

\begin{theorem}
Let $p$ be an odd prime.  The representation $R_2(\Delta,\sigma)$ for $SL_2(\Z/p^2 \Z)$ decomposes into three irreducible 
subrepresentations 
$$  R_2(\sigma) \cong 1\oplus R'_2(\chi_1, \sigma)  \oplus R_2(\chi_{-1}, \sigma)  = \C[\Z/p^2 \Z\times 0].$$
The elements 
\medskip

$$\left\{\sum_{i = 0}^{p-1} (ip,0)\right\} \bigcup \left\{ (j,0)+(p^2-j,0): 1\leq j\leq \frac{p^2-1}2 \right\}$$
\medskip

\noindent are a basis for $R_2(1, \sigma)$ with the first element corresponding to the trivial representation. The basis elements  for  $R_2(-1, \sigma)$ are of two kinds: 

\medskip

\begin{enumerate}
\item  $\left\{(j,0)-(p^2-j,0): (j,p)=1\ \textrm{and}\ 1\leq j\leq \dfrac{p^2-1}2\right\}$ 
\item $\left\{ 2(0,0)-(i,0)-(p^2-i):(i,p)=0\ \textrm{and}\ 1\leq i\leq \dfrac{p^2-1}2 \right\}.$
\end{enumerate}

\medskip

We are throughout identifying elements of $\Z/p^2 \Z\times 0$ with their associated basis elements in $\C[\Z/p^2 \Z\times 0]$. 
\end{theorem}

\medskip

With this, the case of $n=2$ is now complete and we can state the following theorem: 

\medskip

\begin{theorem} Let $p$ be an odd prime. The irreducible representations for $SL_2(\Z/p^2 \Z)$ may first be distinguished as those that arise as the projections of irreducible representations of $SL_2(\Z/p \Z)$ (see Eq. (\ref{eq:proj})) and those that do not. Of the latter kind, select $\Delta_1$ and $\Delta_2$ so that $\left(\dfrac{-\Delta_1}p\right)=1$ and $\left(\dfrac{-\Delta_2}p\right)=-1$ and $\sigma_1$ and $\sigma_2$ so that $\left(\dfrac{\sigma_1}p\right)=1$ and $\left(\dfrac{\sigma_2}p\right)=-1$. Then Table 2.2 
below lists all the non-quotient irreducible representations. 
 \end{theorem}
 
 
\begin{center}
\def\arraystretch{2.5}
\begin{table}[!htbp]
\begin{tabular}{| l | c |  c |}
\hline
$R_k(\chi_\ell, \Delta, \sigma)$ & (Inequivalent) Principal $\chi_\ell$ &  Dimension \\ \hline \hline
$R_0(\chi_\ell, \Delta_1, \sigma_1)$& $\{0<\ell<\dfrac{p(p+1)}{2}: (\ell,p)=1\}$ & $p^2-p$\\
\hline
$R_0(\chi_\ell, \Delta_2, \sigma_1)$& $\{0<\ell<\dfrac{p(p-1)}{2}: (\ell,p)=1\}$ &  $p^2+p$\\
\hline
$R_1(\chi_\ell, \Delta_1, \sigma_1)$& $\{0<\ell<p\}$ &  $\dfrac{p^2-1}{2}$\\
\hline
$R_1(\chi_\ell, \Delta_1, \sigma_2)$& $\{0<\ell<p\}$ &  $\dfrac{p^2-1}{2}$\\
\hline
$R_1(\chi_\ell, \Delta_2, \sigma_1)$& $\{0<\ell<p\}$ &  $\dfrac{p^2-1}{2}$\\
\hline
$R_1(\chi_\ell, \Delta_2, \sigma_2)$&$\{0<\ell<p\}$ &  $\dfrac{p^2-1}{2}$\\
\hline
$R_2(\chi_\ell,\Delta_1,\sigma_1)_1$&$ \{1, -1\}$&$\dfrac{p^2-1}{2}$\\
\hline
$R_2(\chi_\ell,\Delta_1,\sigma_2)_1$&$ \{1, -1\}$&$\dfrac{p^2-1}{2}$\\
\hline
\end{tabular}
\caption{The non-quotient irreducible representations of $SL_2(\Z/p^2\Z)$.}
\end{table}\label{tbl:p2reps}
\end{center}

\pagebreak


\section{Example: Representations and character table for $\text{SL$_2$}(\Z/9 \Z) $}

The representations of $\text{SL$_2$}(\Z/9 \Z) $are realized as those that arise as representations of $\text{SL$_2$}(\Z/3 \Z) $ (quotient representations) and those that don't -- the non-quotient representations. For the former, we refer to \cite{MR-duco} where these are worked out in some detail as an example of a different fast Fourier transform algorithm for the groups  $\text{SL$_2$}(\Z/p\Z) $. For the latter, we use the process detailed above. 

The group $\text{SL$_2$}(\Z/9 \Z) $ has order $648$ and $25$ conjugacy classes with representatives:

\begin{equation}
\left\{\begin{array}{llllllll}
\begin{bmatrix}1&0\\0&1 \end{bmatrix}, & 
 \begin{bmatrix}8&0\\0&8 \end{bmatrix},& 
  \begin{bmatrix}4&3\\3&7 \end{bmatrix},&
   \begin{bmatrix}5&6\\6&2 \end{bmatrix},&
    \begin{bmatrix}0&1\\8&0 \end{bmatrix},&
     \begin{bmatrix}8&8\\8&7 \end{bmatrix},& 
      \begin{bmatrix}1&1\\1&2 \end{bmatrix},\\ \\
       \begin{bmatrix}4&0\\6&7 \end{bmatrix}, &
        \begin{bmatrix}5&0\\3&2 \end{bmatrix}, &
         \begin{bmatrix}1&0\\1&1 \end{bmatrix},&
          \begin{bmatrix}8&0\\8&8 \end{bmatrix},&
           \begin{bmatrix}4&3\\7&1 \end{bmatrix},&
            \begin{bmatrix}5&6\\2&8 \end{bmatrix},& 
             \begin{bmatrix}7&6\\4&1 \end{bmatrix},\\ \\
              \begin{bmatrix}2&3\\5&8 \end{bmatrix},&
               \begin{bmatrix}1&0\\2&1 \end{bmatrix},&
                \begin{bmatrix}8&0\\7&8 \end{bmatrix}, &
                \begin{bmatrix}4&3\\2&4 \end{bmatrix},& 
                  \begin{bmatrix}5&6\\7&5 \end{bmatrix},&
                   \begin{bmatrix}7&6\\2&7 \end{bmatrix},& 
                    \begin{bmatrix}2&3\\7&2 \end{bmatrix},\\ \\
                     \begin{bmatrix}1&3\\0&1 \end{bmatrix},&
                      \begin{bmatrix}8&6\\0&8 \end{bmatrix},&
                       \begin{bmatrix}4&6\\3&7 \end{bmatrix},&
                        \begin{bmatrix}5&3\\6&2 \end{bmatrix}\end{array}\right\}
\end{equation}\label{conjclasses}

\medskip
Since $p=3$ we can choose $\Delta'_1=\sigma_1=1$ and $\Delta'_2=\sigma_2=2$. 
Following the procedure outlined in Section 
2 we separate the representations into cases by $k$: 

\noindent {\bf Case $k=0$}. As seen in Table 2.1 when $k=0$ the principal characters depend only on $\left(\dfrac{-\Delta_i}{p}\right)$. When $\Delta=1$ there are eight principal characters, indexed by the equivalence classes
$$ \{1,11\},\{2,10\},\{4,8\},\{5,7\}.$$
Selecting one character from each we get  the set of indices $\{1,2,4,5\}$ which corresponds to four representations of dimension 6. When $\Delta=2$ there are four principal characters, indexed by $\{\{1,5\},\{2,4\}\}$.  Class representatives are given by indices  $\{1,2\}$ producing two representations of dimension 12. 
\noindent {\bf Case $k=1$.} In this case the principal characters are independent of $\Delta$ and $\sigma$. There are four characters $\{\{1,5\},\{2,4\}\}$ (grouped according to equivalence class) so that  $\{1,2\}$ give class representatives. Independently choosing representative character, $\Delta$, and $\sigma$ determines  eight four-dimensional representations. 
\noindent {\bf Case $k=2$.}  In this case we follow the procedure described in Section 2.3.1. This gives a a basis $$\{(1,0)-(8,0),(2,0)-(7,0),(3,0)-(6,0),(4,0)-(5,0)\}$$ for the non-trivial part of $R_2(\sigma_j,\chi_1)$ and a basis $$\{(1,0)+(8,0),(2,0)+(7,0),2(0,0)-(3,0)-(6,0),(4,0)+(5,0)\}$$ for $R_2(\sigma_j,\chi_{-1})$ to obtain the final four representations of dimensions four.

In total there are 18 non-quotient representations which along  with the seven irreducible representations of $SL_2(\Z/3\Z)$ gives the  25 irreducible representations of $SL_2(\Z/3^2 \Z)$ (recall that the number of irreducible representations is equal to the number of conjugacy classes). These representations are summarized in Table \ref{tbl:Example} below. The character table for the new representations is also constructed below with the rows indexed by conjugacy classes and the columns indexed by representations. 
\begin{center}
\def\arraystretch{2.5}\label{tbl:Example}
\begin{table}[!htbp]
\begin{tabular}{c}
{

\setcounter{MaxMatrixCols}{25}
\resizebox{\linewidth}{!}{
$         	
\begin{bmatrix}  
R_0(1, 1, 1) & R_0(2, 1, 1) & R_0(4, 1, 1) & R_0(5, 1, 1) & R_0(1, 2, 1) & R_0(2, 2, 1) & R_1(1, 1, 1) & R_1(2, 1, 1) & R_1(1, 2, 1) & R_1(2, 2, 1) & R_1(1, 1, 2) & R_1(2, 1, 2) & R_1(1, 2, 2) & R_1(2, 2, 2) &R_2(0, 1) &R_2(1, 1) & R_2(0,2) & R_2(1,2) \\ 

6							& 6		& 6			& 6							& 12		& 12		& 4						& 4							& 4						& 4						&4						& 4						& 4						& 4						& 4						& 4						& 4						& 4						\\
-6							& 6		& 6			& -6							& -12	& 12		& 4						& -4							& 4						& -4						&4						& -4						& 4						& -4						& 4						& -4						& 4						& -4						\\
3							& 3		& 3			& 3							& 0		& 0		& -2						& -2							& -2						& -2						&-2						& -2						& -2						& -2						& -2						& -2						& -2						& -2						\\
-3							& 3		& 3			& -3							& 0		& 0		& -2						& 2							& -2						& 2						&-2						& 2						& -2						& 2						& -2						& 2						& -2						& 2						\\
0							& 0		& 0			& 0							& -3		& -3		& 1						& 1							& 1						& 1						&1                                        	& 1						& 1						& 1						& 1						& 1						& 1						& 1						\\
0							& 0		& 0			& 0							& 3		& -3		& 1						& -1							& 1						& -1						&1                                        	& -1						& 1						& -1						& 1						& -1						& 1						& -1						\\
-3							& -3		& -3			& -3							& 3		& 3		& -\zeta_{3} + 2\zeta_{3}^2	& -\zeta_{3} + 2\zeta_{3}^2		& 2\zeta_{3} - \zeta_{3}^2		& 2\zeta_{3} - \zeta_{3}^2		&-\zeta_{3} + 2\zeta_{3}^2    	& -\zeta_{3} + 2\zeta_{3}^2	& 2\zeta_{3} - \zeta_{3}^2		& 2\zeta_{3} - \zeta_{3}^2		& -\zeta_{3} + 2\zeta_{3}^2	& -\zeta_{3} + 2\zeta_{3}^2	& 2\zeta_{3} - \zeta_{3}^2		& 2\zeta_{3} - \zeta_{3}^2		\\
3							& -3		& -3			& 3							& -3		& 3		& -\zeta_{3} + 2\zeta_{3}^2	& \zeta_{3} - 2\zeta_{3}^2			& 2\zeta_{3} - \zeta_{3}^2		& -2\zeta_{3} + \zeta_{3}^2	&-\zeta_{3} + 2\zeta_{3}^2     	& \zeta_{3} - 2\zeta_{3}^2 		& 2\zeta_{3} - \zeta_{3}^2		& -2\zeta_{3} + \zeta_{3}^2	& -\zeta_{3} + 2\zeta_{3}^2	& \zeta_{3} - 2\zeta_{3}^2		& 2\zeta_{3} - \zeta_{3}^2		& -2\zeta_{3} + \zeta_{3}^2)	\\
-3							& -3		& -3			& -3							& 3		& 3		& 2\zeta_{3} - \zeta_{3}^2		& 2\zeta_{3} - \zeta_{3}^2			& -\zeta_{3} + 2\zeta_{3}^2	& -\zeta_{3} + 2\zeta_{3}^2	&2\zeta_{3} - \zeta_{3}^2   	& 2\zeta_{3} - \zeta_{3}^2		& -\zeta_{3} + 2\zeta_{3}^2	& -\zeta_{3} + 2\zeta_{3}^2	& 2\zeta_{3} - \zeta_{3}^2		& 2\zeta_{3} - \zeta_{3}^2		& -\zeta_{3} + 2\zeta_{3}^2	& -\zeta_{3} + 2\zeta_{3}^2	\\
3							& -3		& -3			& 3							& -3		& 3		& 2\zeta_{3} - \zeta_{3}^2		& -2\zeta_{3} + \zeta_{3}^2		& -\zeta_{3} + 2\zeta_{3}^2	& \zeta_{3} - 2\zeta_{3}^2		&2\zeta_{3} - \zeta_{3}^2   	& -2\zeta_{3} + \zeta_{3}^2 	& -\zeta_{3} + 2\zeta_{3}^2	& \zeta_{3} - 2\zeta_{3}^2		& 2\zeta_{3} - \zeta_{3}^2		& -2\zeta_{3} + \zeta_{3}^2	& -\zeta_{3} + 2\zeta_{3}^2	& \zeta_{3} - 2\zeta_{3}^2		\\
0							& -2		& 2			& 0							& 0		& 0		& 0						& 0							& 0						& 0						&0						& 0						& 0						& 0						& 0						& 0						& 0						& 0						\\
-\zeta_{12}^7 + \zeta_{12}^{11}		& 1		& -1			& \zeta_{12}^7 - \zeta_{12}^{11}	& 0		& 0		& 0						& 0							& 0						& 0						&0						& 0						& 0						& 0						& 0						& 0						& 0						& 0						\\
\zeta_{12}^7 - \zeta_{12}^{11}		& 1		& -1			& -\zeta_{12}^7 + \zeta_{12}^{11}	& 0		& 0		& 0						& 0							& 0						& 0						&0						& 0						& 0						& 0						& 0						& 0						& 0						& 0						\\
0							& 0		& 0			& 0							& 0		& 0		& -\zeta_{3}				& -\zeta_{3}					& -\zeta_{3}^2				& -\zeta_{3}^2				&-\zeta_{3}^2				& -\zeta_{3}^2				& -\zeta_{3}				& -\zeta_{3}				& 1						& 1						& 1						& 1						\\
0							& 0		& 0			& 0							& 0		& 0		& -\zeta_{3}				& \zeta_{3}					& -\zeta_{3}^2				& \zeta_{3}^2				&-\zeta_{3}^2				& \zeta_{3}^2				& -\zeta_{3}				& \zeta_{3}				& 1						& -1						& 1						& -1						\\
0							& 0		& 0			& 0							& 0		& 0		& -1						& -1							& -1						& -1						&-\zeta_{3}				& -\zeta_{3}				& -\zeta_{3}^2				& -\zeta_{3}^2				& \zeta_{3}^2				& \zeta_{3}^2				& \zeta_{3}				& \zeta_{3}				\\
0							& 0		& 0			& 0							& 0		& 0		& -1						& 1							& -1						& 1						&-\zeta_{3}				& \zeta_{3}				& -\zeta_{3}^2				& \zeta_{3}^2				& \zeta_{3}^2				& -\zeta_{3}^2				& \zeta_{3}				& -\zeta_{3}				\\
0							& 0		& 0			& 0							& 0		& 0		& -\zeta_{3}^2				& -\zeta_{3}^2					& -\zeta_{3}				& -\zeta_{3}				&-1						& -1						& -1						& -1						& \zeta_{3}				& \zeta_{3}				& \zeta_{3}^2				& \zeta_{3}^2				\\
0							& 0		& 0			& 0							& 0		& 0		& -\zeta_{3}^2				& \zeta_{3}^2					& -\zeta_{3}				& \zeta_{3}				&-1						& 1						& -1						& 1						& \zeta_{3}				& -\zeta_{3}				& \zeta_{3}^2				& -\zeta_{3}^2				\\
0							& 0		& 0			& 0							& 0		& 0		& -\zeta_{3}^2				& -\zeta_{3}^2					& -\zeta_{3}				& -\zeta_{3}				&-\zeta_{3}				& -\zeta_{3}				& -\zeta_{3}^2				& -\zeta_{3}^2				& 1						& 1						& 1						& 1						\\
0							& 0		& 0			& 0							& 0		& 0		& -\zeta_{3}^2				& \zeta_{3}^2					& -\zeta_{3}				& \zeta_{3}				&-\zeta_{3}				& \zeta_{3}				& -\zeta_{3}^2				& \zeta_{3}^2				& 1						& -1						& 1						& -1						\\
0							& 0		& 0			& 0							& 0		& 0		& -\zeta_{3}				& -\zeta_{3}					& -\zeta_{3}^2				& -\zeta_{3}^2				&-1						& -1						& -1						& -1						& \zeta_{3}^2				& \zeta_{3}^2				& \zeta_{3}				& \zeta_{3}				\\
0							& 0		& 0			& 0							& 0		& 0		& -\zeta_{3}				& \zeta_{3}					& -\zeta_{3}^2				& \zeta_{3}^2				&-1						& 1						& -1						& 1						& \zeta_{3}^2				& -\zeta_{3}^2				& \zeta_{3}				& -\zeta_{3}				\\
0							& 0		& 0			& 0							& 0		& 0		& -1						& -1							& -1						& -1						&-\zeta_{3}^2				& -\zeta_{3}^2				& -\zeta_{3}				& -\zeta_{3}				& \zeta_{3}				& \zeta_{3}				& \zeta_{3}^2				& \zeta_{3}^2				\\
0							& 0		& 0			& 0							& 0		& 0		& -1						& 1							& -1						& 1						&-\zeta_{3}^2				& \zeta_{3}^2				& -\zeta_{3}				& \zeta_{3}				& \zeta_{3}				& -\zeta_{3}				& \zeta_{3}^2				& -\zeta_{3}^2				\\
 \end{bmatrix}
$
}
}\\
$\;$
\end{tabular}\caption{Character Table for the non-quotient irreducible representations of $SL_2(\Z/9\Z)$ The rows are indexed according to the conjugacy classes as listed in Eq. (\ref{conjclasses}) and the columns correspond to the characters as listed previously according to $k= 0,1,2$.}

\end{table}
\end{center}

\begin{center}
\def\arraystretch{2.5}
\begin{table}[!htbp]
\begin{tabular}{| l | c | c  || l | c | c|}
\hline
\multicolumn{3}{|c||}{Quotient Representations}&\multicolumn{3}{|c|}{Non-Quotient Representations}\\
\hline
Representation & Type &   Dimension  &$R_k(\chi, \Delta, \sigma)$ & Principal $\chi_\ell$ &   Dimension \\
\hline \hline                                                        
$\rho_0$ &    PS &  $1$ &    $R_0(\chi, 2, 1)$ & $\{1,2\}$ &  12 \\
\hline                               
$\rho_{012}$ & PS  & 3 & $R_0(\chi, 1, 1)$ & $\{1,2,4,5\}$ & 6\\
\hline                                
$\rho_{01}$&   PS &  $2$&       $R_1(\chi, 1, 1)$&  $\{1,2\}$ &   4\\                      
\hline                       
$\rho_{02}$&  PS   &   $2$&      $R_1(\chi, 1, 2)$&  $\{1,2\}$ &   4\\                                       
\hline                                                                            
$\rho_{12}$&  DS   &   $2$&      $R_1(\chi, 2, 1)$&  $\{1,2\}$ &   4\\                                    
\hline                                                                            
$\rho_{1}$&   DS  &   $1$&          $R_1(\chi, 2, 2)$&  $\{1,2\}$ &   4\\                              
\hline                                                                            
$\rho_{2}$&   DS &   $1$&               $R_2(\chi, 1, 1)$&  $\{1,2\}$ &   4\\                            
\hline                                                                                  
           &   &      &       $R_2(\chi, 1, 2)$&  $\{1,2\}$ &   4\\   
\hline
\end{tabular}\label{tbl:z9reps}           
\caption{The  irreducible representations of $SL_2(\Z/3^2 \Z)$. The ``Type" for the quotient representations are either {\em principal series} (PS) or {\em discrete series} (DS). The notation for quotient representations as in \cite{MR-duco}, with details for the construction in \cite{LR92}. The notation for the non-quotient representations is as in Section 2.}
\end{table}
\end{center}

\pagebreak

%
%

\pagebreak
\clearpage
%
%
%
%
%
%
%
%
%
%
%
%
%
\section{Numerical Experiments}

The ability to construct representations of $\slpsq$ enables some experimentation on the spectra of Cayley graphs on these groups. Recall that for any group $G$ and generating set $S = \{s_1,\dots, s_k\}$ such that $S=S^{-1}$, the Cayley graph $X(G,S)$ is the network with node set $G$ such that any $g\in G$ is connected to the each of the elements $gs_i$. The {\em spectrum} of $X(G,S)$ is the set of eigenvalues of the related {\em adjacency matrix} $A(G,S)$, the $|G|\times |G|$ matrix with 
$a_{g,gs_i} = 1$ for $i=1,\dots,k$ and $0$ otherwise. 

Note that 
$$A(G,S) = \sum_{s\in S} \rho(s)$$
where $\rho$ is the right multiplication representation of $G$. By Wedderburn's Theorem, 
$$\mbox{spectrum}\left[A(G,S)\right] = \bigcup_{\eta \in \widehat{G}} \mbox{spectrum}\left[\widehat{\delta_S}(\eta)\right]$$
where the sum is over a complete set of irreducible representations of $G$ and $\widehat{\delta_S}(\eta)$ denotes the {\em Fourier transform} of $\delta_S$, the characteristic function for the set $S$, at the irreducible representation $\eta$ of $G$, defined as 
$$\widehat{\delta_S}(\eta) = \sum_{s\in S}\eta(s).$$
Note that the {\em normalized} form would compute ${\frac{1}{|S|}} \widehat{\delta_S}(\eta)$. 

Since any irreducible representation can have dimension at most $|G|^{1/2}$, if $|G|$ is large, this observation can enable eigenvalue calculations for $X(G,S)$ that might have been intractable if attempted directly using $A(G,S)$ (see e.g., \cite{LR92,LR93}). In particular, for $G=\slpsq$ this is the difference between finding the eigenvalues of a matrix of order $p^6$ versus roughly $p^2$ matrices of order $p^2$.

Similarly, suppose $G$ generated by $S$ is acting on a set $\Gamma$. The Cayley graph on $\Gamma$ has point set $\Gamma$ and connects $\gamma$ to $s\gamma $ for any $s\in S$. If $\Gamma = G/B$ for a subgroup $B$ then the spectrum of the Cayley graph is equal to the union of the spectra over any representations that appear in the induction of the trivial representation from $B$ to $G$. 

\subsection{Monochromatic Eigenvalues} Computations and conjectures concerning the expansion properties and spacings of eigenvalues of Cayley graphs on $\slp$ have been well studied for many years.  Recently, there has been renewed interest in this problem as a result of connections to quantum chaos. See the introduction of \cite{RS} for a discussion of the progression of these results for $\slp$ and interpretations of earlier results in terms of modern terminology. 

In particular, in this setting, a process of desymmetrization is applied to the spectrum by coloring each eigenvalue with the irreducible representation to which it belongs. Eigenvalues that only arise in the spectrum of a single representation are called {\bf monochromatic} and these are the values whose behavior is of statistical interest. For our experiments below, this leads us to consider the eigenvalues of irreducible representations individually rather than collectively. 

\subsection{Spectra of some Cayley graphs on $SL_2(\Z /p^2 \Z)$}
Here we consider Cayley graphs for three families of generators, defined below, which have been previously used for similar experiments on $SL_2(\Z/p\Z)$. The first two were originally introduced in \cite{LR92} and the ``Selberg family" in \cite{LR99,RS}. We expect based on the previous results that $G_1$ and $G_2$ should be associated to non--expanders while $G_3$ should give second eigenvalues much closer to the Ramanujan bound.  

%
%


$$\begin{array}{rcccc}

G_1 = &  \left\{\begin{bmatrix}
1&1\\0&1
\end{bmatrix},\begin{bmatrix}
1&-1\\0&1
\end{bmatrix} ,\begin{bmatrix}
0&1\\-1&0
\end{bmatrix} ,\begin{bmatrix}
0&-1\\1&0
\end{bmatrix}\right\} \\ \\
G_2 = & \left\{\begin{bmatrix}
1&1\\-1&0
\end{bmatrix},\begin{bmatrix}
0&-1\\1&1
\end{bmatrix} ,\begin{bmatrix}
0&1\\-1&0
\end{bmatrix} ,\begin{bmatrix}
0&-1\\1&0
\end{bmatrix}\right\}  \\ \\
G_3 =& \left\{\begin{bmatrix}
1&2\\0&1
\end{bmatrix},\begin{bmatrix}
1&-2\\0&1
\end{bmatrix} ,\begin{bmatrix}
1&0\\-2&1
\end{bmatrix} ,\begin{bmatrix}
1&0\\2&1
\end{bmatrix} \right\}\\

%
\end{array}$$

We begin by restricting our attention to the particular case of $SL_2(\Z /25 \Z)$ to provide a more detailed analysis of the spectra. This group has 15,000 elements and 40 non--quotient irreducible representations. The parameterization and ordering of these representations along the y--axis in Figures 1--3 is given in Appendix A. 

For each of the non--quotient irreducibles we begin by extracting the leading eigenvalues to get a sense of the properties of these values over the parameter choices. The clearest behavioral pattern is determined by the parity of the associated character. Note that the adjacent elements (from bottom to top) that appear to have the same magnitude are due to skipping multiples of $5$ in the choice of characters for the $k=0$ case. In many of the examples for small primes that we computed, including ones not displayed here, the largest eigenvalue across the irreducibles occurred in the $k=2$ representation where $\left(\frac{\sigma}{p}\right)=\chi=1$.

\begin{figure}[!h]
\subfloat[$G_1$]{\includegraphics[height=2.25in]{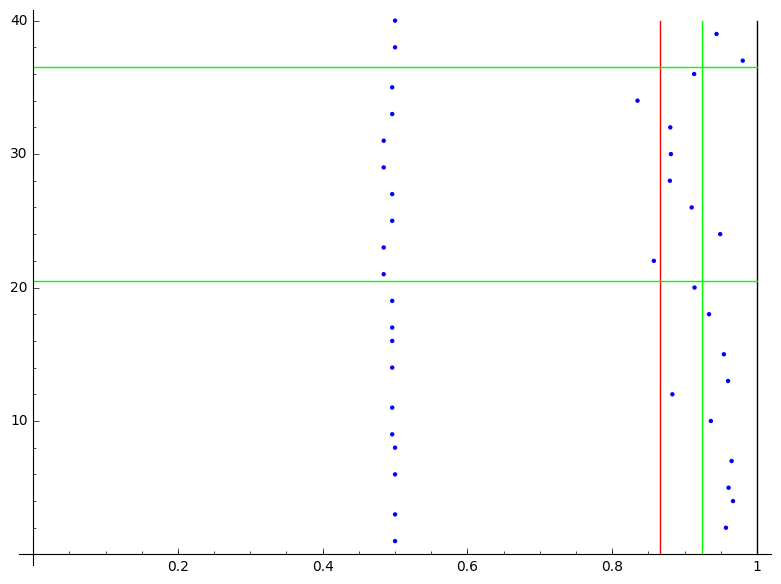}} \\
\subfloat[$G_2$]{\includegraphics[height=2.25in]{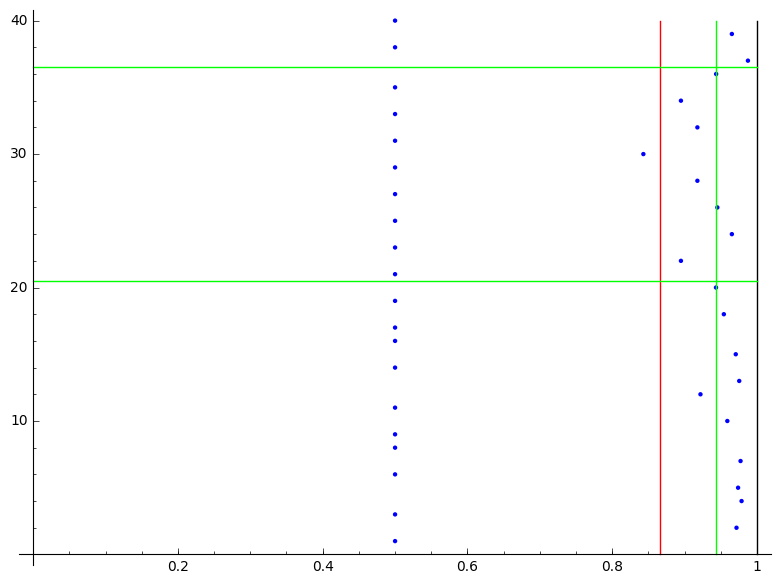}} 
\subfloat[$G_3$]{\includegraphics[height=2.25in]{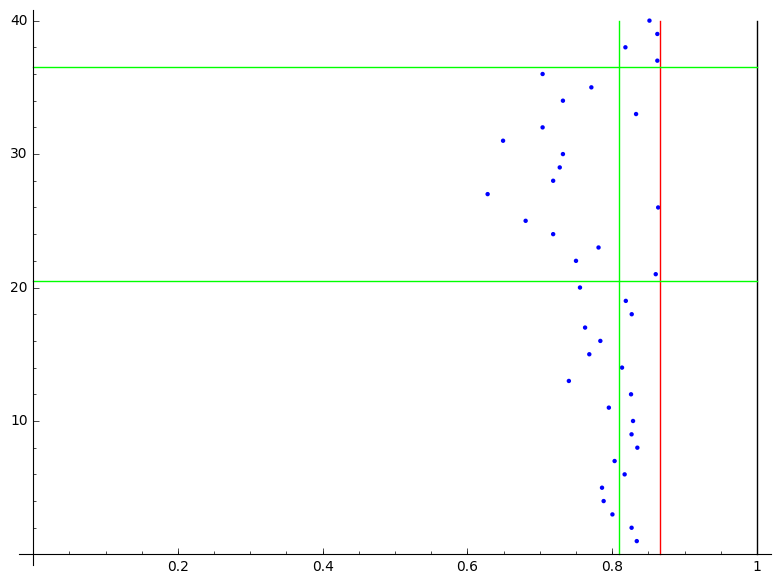}}
\caption{Plots of the second largest eigenvalues of the non--quotient irreducible representations for $SL_2(Z/25Z)$ ordered as in the previous figure.  The green lines separate representations for different values of $k$ while the vertical green, red, and black lines marks the $\lambda_1$ of $SL_2(\Z /5\Z)$ for the generating set, the normalized Ramanujuan bound, and leading eigenvalue, respectively. 
 }
\end{figure}

 
The full spectrum for each of the non--quotient irreducible representations for the generating sets $G_1$, $G_2$, and $G_3$ are presented below. These plots begin to give a sense of the chromatic structure of the eigenvalues and suggest that for these representations it may be possible to describe the genericity in terms of $k$ and the parity of $\chi$. The structure suggested by the eigenvalues that appear in all or most of the representations such as 1 for generating set $G_2$ presents an avenue for future study. 
 
 \vfill
\begin{figure}[!h]
\subfloat[$G_1$]{\includegraphics[height=2.5in]{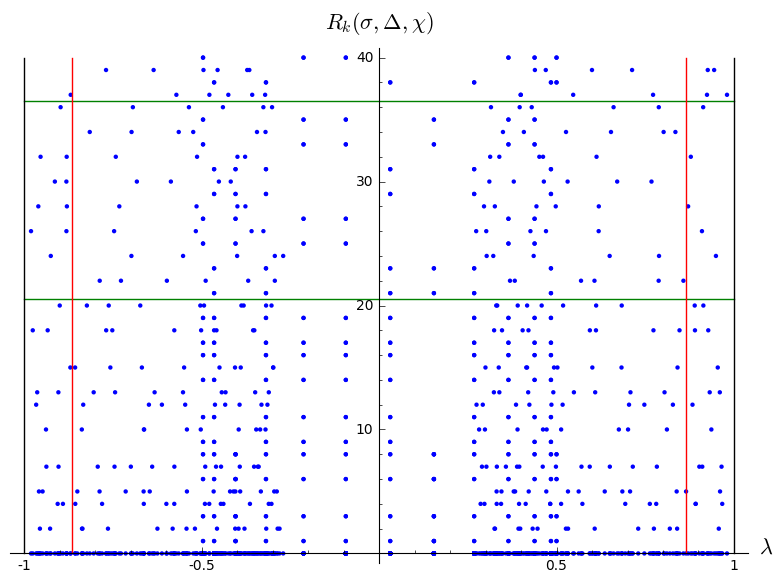}} \\
\subfloat[$G_2$]{\includegraphics[height=2.5in]{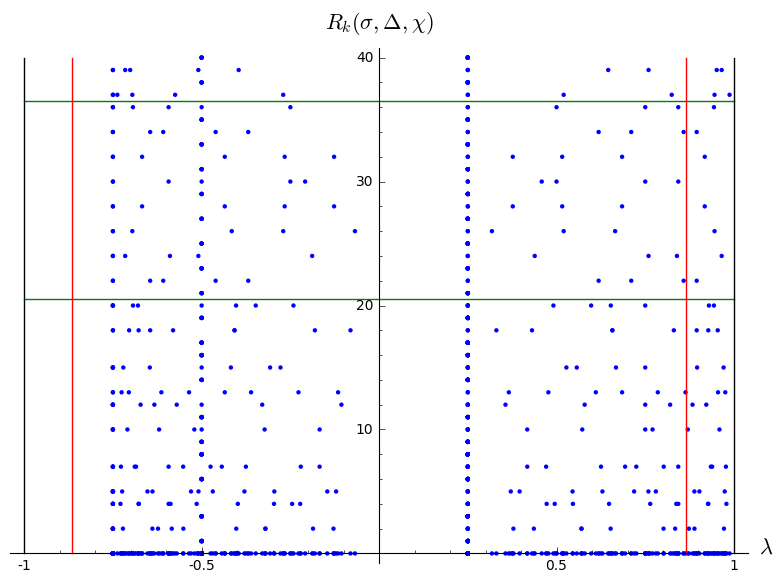}} 
\subfloat[$G_3$]{\includegraphics[height=2.5in]{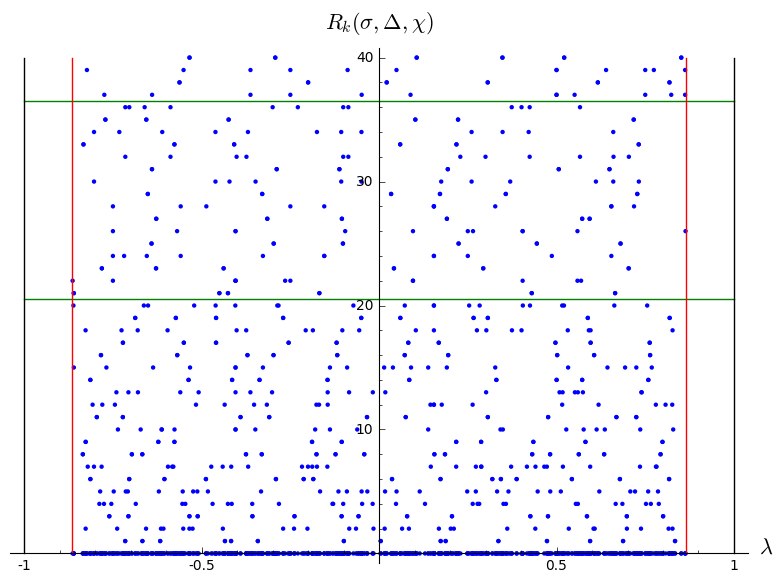}}
\caption{Plots of the non--quotient eigenvalues of the full Cayley graphs for $SL_2(Z/25Z)$ for the generating sets described above, separated by irreducible representation. The values at height $0$ represent the full spectrum for the non--quotient representations while the horizontal lines above each correspond to a single irreducible. The ordering along the y--axis is described in Appendix A. The green lines separate representations for different values of $k$ while the vertical black and red lines mark the leading and Ramanajuan values. }
\end{figure}

\vfill


Next, we compute some of the Cayley spectra for larger $p$ and $n$, highlighting some connections to previous work and directions for future research. For the groups $SL_2(\Z /p^2 \Z)$ with $p\in\{3,5,7,11\}$ we compute the eigenvalues corresponding to the $k=2$ spaces, which have the connection to the leading eigenvalues for the entire representation. These should be compared to the top four representations in Figures 1 and 2, combining the two choices each for of $\sigma$ and $\chi$. Notice that for the generating sets $G_1$ and $G_2$ the leading eigenvalue of the Cayley graph is contained in this parameter space, and that these graphs are not expanders for any $p$ we computed.  while for $G_3$ we observe the Ramanujuan gap. We conjecture that this is the case for all $p$, identifying these as an important family of representations. Following \cite{RS} we also conjecture that $G_3$ should remain an expander as $p$ grows. 

\begin{figure}[!h]
\subfloat[$G_1$]{\includegraphics[height=2.5in]{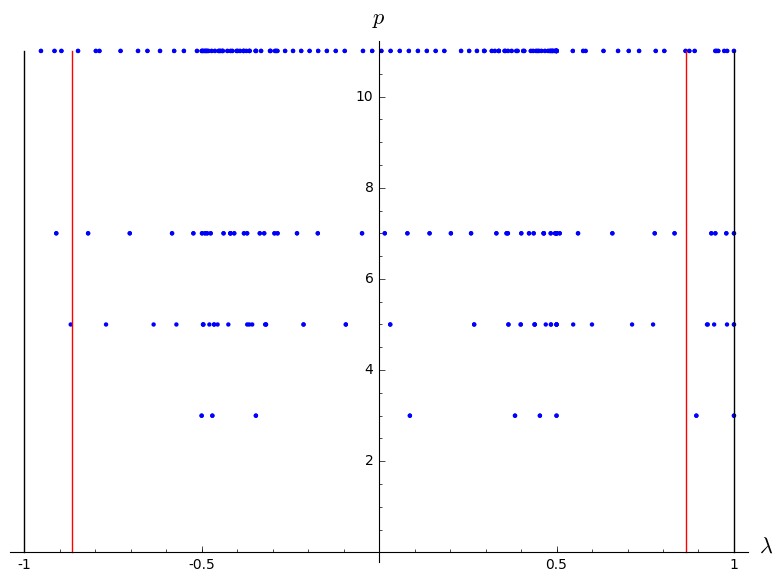}}\\
\subfloat[$G_2$]{\includegraphics[height=2.5in]{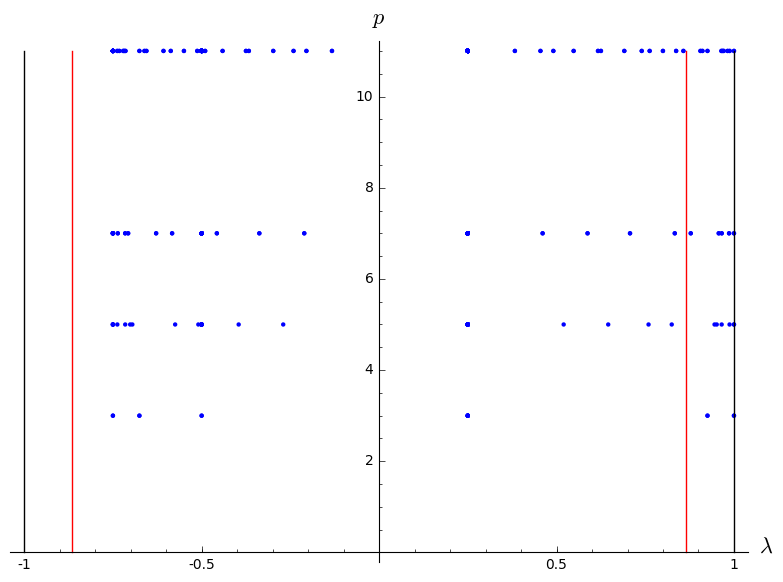}}
\subfloat[$G_3$]{\includegraphics[height=2.5in]{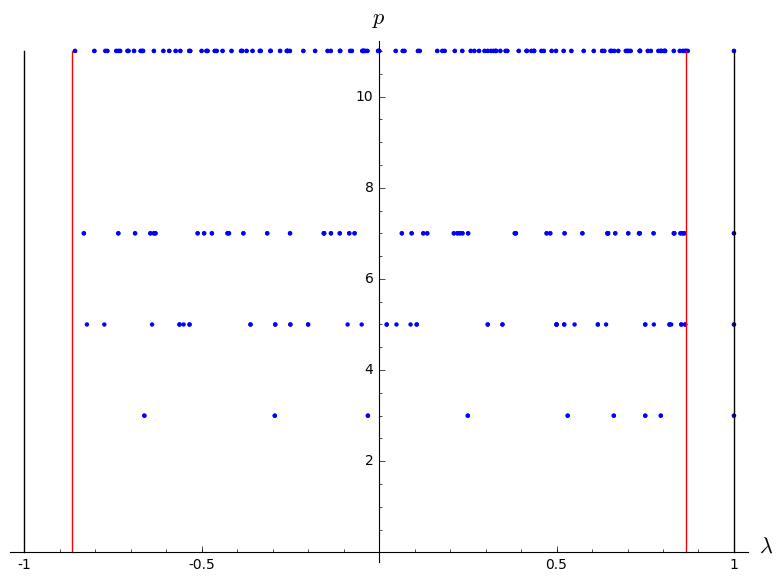}}

\caption{Cayley eigenvalues of the $k=2$ spaces for  $SL_2(\Z /p^2 \Z)$ with $p\in\{3,5,7,11\}$ for each of our generating sets. The y--axis represents the $p$ value and the vertical red line denotes the maximal (normalized) eigenvalue of $1$, which also occurs in the full $R_2(\sigma)$ space as a copy of the identity cf. Section 2.3.  }
\end{figure}

%
%
%
%

Finally, we construct some of the Cayley spectra for $SL_2(\Z /125 \Z)$. The question of how the eigenvalues in the various parameter spaces nest moving from $n=2$ to $n=3$, is interesting.  Particularly in terms of identifying new monochromatic eigenvalues, since they can arise in two ways, as either   a monochromatic value for $n=2$ that is not duplicated when lifting or an entirely new value that occurs at $n=3$. 

This figure displays the eigenvalues for the $\Delta=\sigma=1$ and $\chi\in\{1,2\}$ spaces for $k=1$ and $k=2$. Notice that for $G_2$ the $k=\chi=1$ case matches exactly the two eigenvalues shown in Figure 3(B), while the behavior is much less clear for some of the other representations. Determining natural correspondences between these spaces is an interesting question for future work but the structure here leads us to conjecture that the parameterization of irreducibles given by the $R_k(\sigma, \Delta, \chi)$ plays a similar role both as $n$ and $p$ increase.

\begin{figure}[!h]

\subfloat[$G_1$]{\includegraphics[height=2.5in]{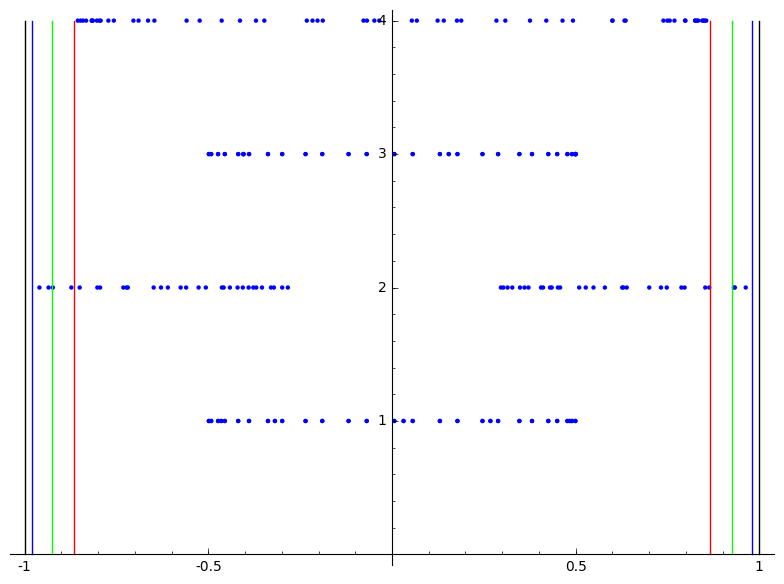}}\\
\subfloat[$G_2$]{\includegraphics[height=2.5in]{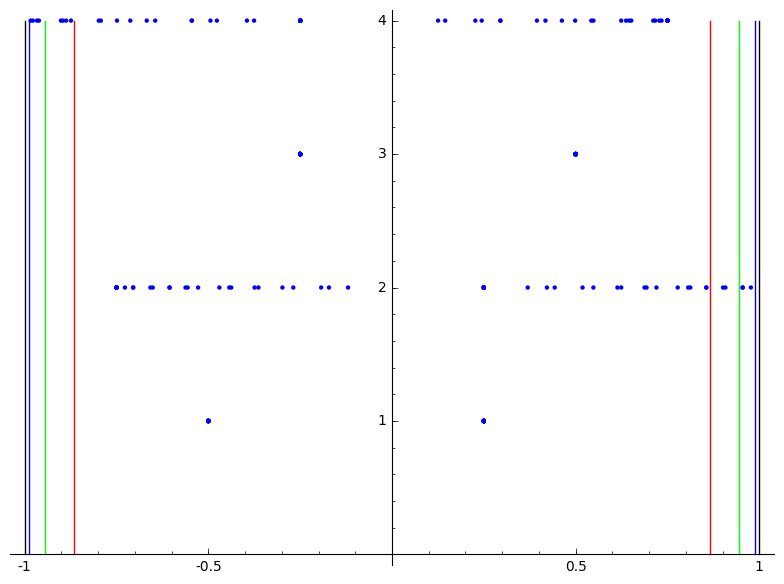}}
\subfloat[$G_3$]{\includegraphics[height=2.5in]{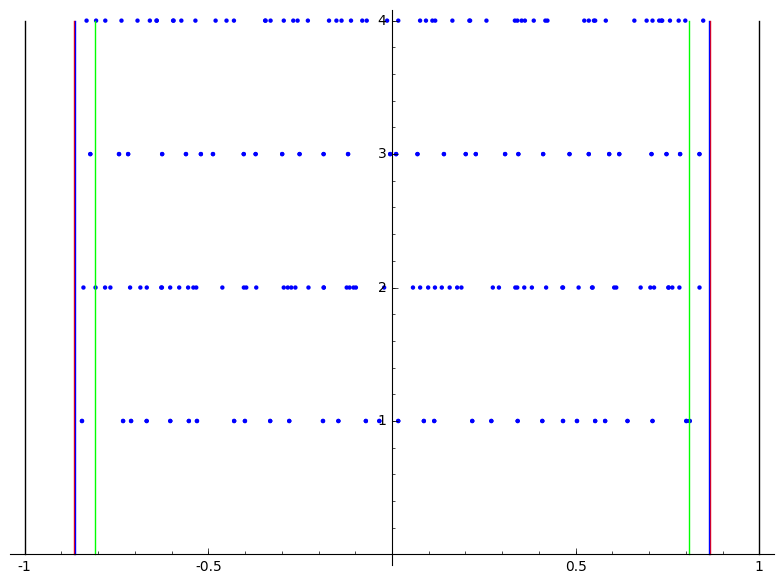}}
\caption{Cayley eigenvalues for the $k=1$ and $k=2$ spaces for   $SL_2(\Z /125 \Z)$ with $\Delta=\sigma=1$ and $\chi\in\{1,2\}$ for each space. The $\lambda_1$ values for $SL_2(\Z /25 \Z)$ (blue) and $SL_2(\Z /5 \Z)$ (green) as well as the Ramanujuan bound (red) are marked with vertical lines.  Parameters for the y--axis are given in Appendix A. }
\end{figure}

The table below compares the $\lambda_1$ values for $p=5$ and $k=1,2,3$ for the three generating sets we consider. 

\begin{table}[!h]
\begin{tabular}{|l|c|c|c|}
\hline
$n$&1&2&3\\
\hline
\hline
$G_1$&3.699&3.922&$>3.852$\\
\hline
$G_2$&3.774&3.950&$>3.935$\\
\hline
$G_3$&3.236&3.454&$>3.389$\\
\hline
\end{tabular}
\caption{Comparison of $\lambda_1$ for Cayley graphs on $SL_2(\mathbb{Z}/p^n\mathbb{Z})$. }
\end{table}

\subsection{Spectrum of the graph on the ``projective line" $G/B$} Motivated by \cite{RS} we   considered the ``projective'' Schreier graph formed by the action of the generators on the cosets of $G/B$, where $B$ is the upper triangular Borel subgroup. 
We first computed the second largest eigenvalue for the generating sets described in the previous section for primes less than 50 Figure~5 (a). As observed in \cite{LR92}, the eigenvalues appear to converge to a limit quickly. Note that these provide a lower bound for the second largest eigenvalue of the full Cayley graph with the same generators and hence $G_1$ and $G_2$ seem to share the behavior shown in Figures 4, 5, and 6 of \cite{LR92}. The case of $G_3$  for $SL_2(Z/pZ)$ was considered in Figure 2 of \cite{RS}  and our result appears similar, although over a much smaller range of primes. For $G_1$ the leading eigenvalue of the Schreier graph was equal to the leading eigenvalue of the full Cayley graph. We conjecture that this continues for larger $p$. 

We next generated 100 pairs of generating elements $\{s,t\}$ for $SL_2(Z/p^2Z)$ for $3\leq p\leq 19$ and computed the second largest eigenvalues for the Cayley graphs generated by $\{s,s^{-1},t,t^{-1}\}$. These are displayed in Figure~1(b). This corresponds to the experiments considered in Figure 4 of \cite{RS} and Figure 8 of \cite{LR92}. The mean of the values appears to be approaching the (normalized) Ramanujuan bound of $2\sqrt{3}/4$ on average as $p$ increases which would correspond to the results in \cite{RS} for $\slpn$, thus we expect the same behavior to continue in this new setting. 

\begin{figure}[!h] 
\label{fig:GmodB}
\subfloat[Fixed Generators]{\includegraphics[height=2in]{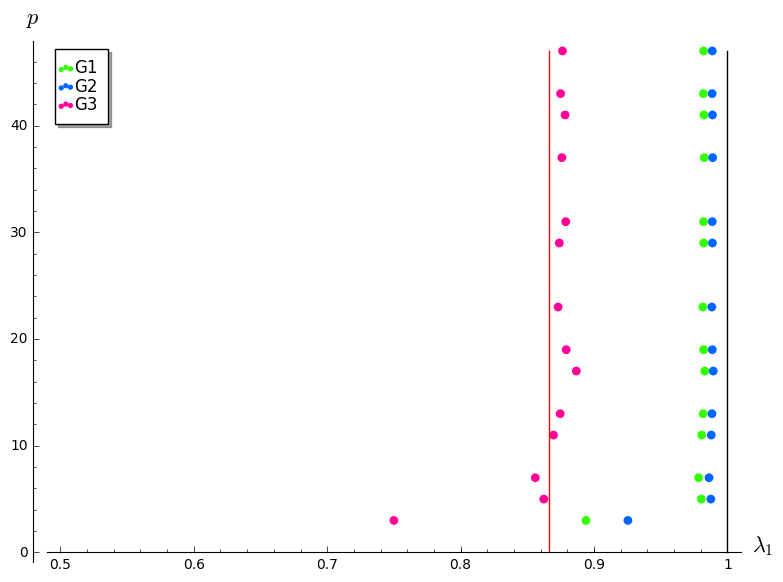}}\qquad
\subfloat[Random Generators]{\includegraphics[height=2in]{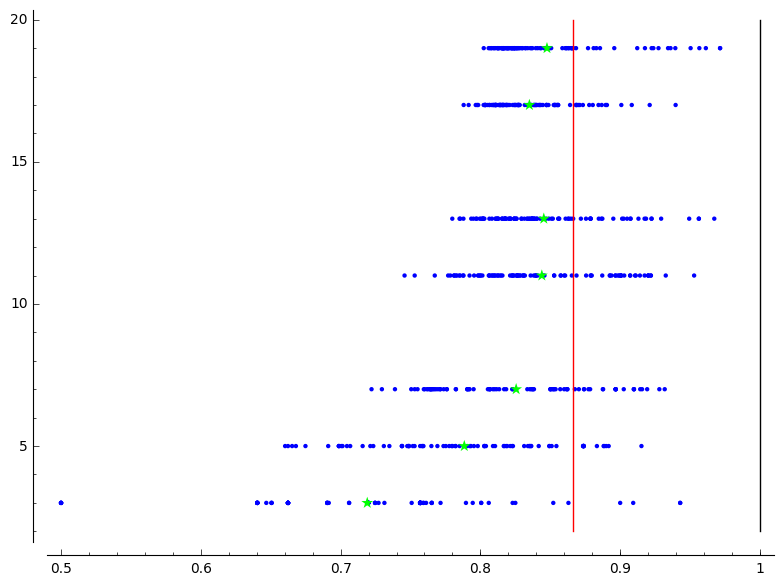}}
\caption{Eigenvalues associated to the Schreier graph $SL_2(Z/p^2Z)/B$. Figure (a) shows the behavior of the second largest eigenvalues for the three generating sets given in Section 4.1 as $p$ varies, normalized to lie in $[-1,1]$. The red lines represent the largest normalized eigenvalue $1$ and the corresponding Ramanujuan bound $\frac{\sqrt{3}}{2}$. Of particular interest is the behavior of $G_3$ (pink dots) which also generates nearly Ramanajuan graphs in the $SL_2(Z/pZ)/B$ case \cite{RS}. Figure (b) shows plots of the second largest eigenvalues for 100 random generating pairs, also  normalized to lie in $[-1,1]$. The red lines represent the largest eigenvalue $4$ and the Ramanujuan bound. For $p=19$ the mean of the (un--normalized) values is approximately 3.41 and it is natural to conjecture that as $p$ increases the mean converges to $2\sqrt{3}$. The mean for each $p$ is marked with a green $*$.   }
\end{figure}


\section{Directions for Future Work}
The construction and experiments presented here provide a baseline for future work. Computationally, more efficient methods are needed to handle the large sizes of the representation spaces as $n$ and $p$ grow. Generating tractable bases for $k=n$ when $k>2$ and canonical choices of characters for the case when $C$ is non--cyclic are also interesting problems. With respect to the representations themselves, determining which parameter sets have related, nested, or preserved eigenvalues as $n$ grows for fixed $p$ is important, as well as determining which families of monochromatic eigenvalues have the same expected values.


Additionally, our experiments provide support for generalizing the results developed for $SL_2(\mathbb{Z}/p\mathbb{Z})$ in \cite{RS} to $SL_2(\mathbb{Z}/p^n\mathbb{Z})$. In particular, we conjecture that for large $p$, Cayley graphs on random generating sets for $SL_2(\mathbb{Z}/p^n\mathbb{Z})$ should be almost--Ramanajuan and that the same should hold for the for the analogue of the projective line. Similarly, we conjecture that the monochromatic eigenvalues should be independent and be centered around the Ramajuan bound  as $p$ increases. 




\section*{Acknowledgements} Thanks to Alex Lubotzky, Peter Sarnak, Justin Troyka, Sam Schiavone, Zachary Garvey, Thomas Shemanske and John Voight for helpful conversations and encouragements. 

\bibliographystyle{plain}
\bibliography{SL2p2}


\appendix

\section{Additional Mathematical Context}
The explicit construction presented in Section 2 obscures some of the motivating formalism that lead to the discovery of these methods. In this appendix we provide some additional context for the construction, using the same notation as in Section 2. 
In the literature, the representations of $\slpn$ are defined on the vector space of complex functions $f : G \to \C$. We work with $\C G$ simply to produce explicit matrices. This space $\C G$ is dual to this space of complex functions on $G$ which causes our statements to be reversed from those in the standard literature.

In addition, the construction in \cite{nobs1976irreduziblen}, \cite{nobs1976irreduziblenp1} works with an underlying quadratic form described below instead of defining the multiplication directly on $G$. In order to construct this form, choose $\Delta'$ and $\sigma$ to be two integers such that $(\Delta',p) = (\sigma,p) =1$. Let $\Delta = p^k \Delta' $ and define a quadratic form $Q$:

\[ \begin{aligned} 
Q : M & \to \Z/p^n\Z   &  \\
(x,y) & \mapsto x^2 + \Delta y^2   &   
\end{aligned} \]

\noindent and its associated bilinear form $B$

\[ \begin{aligned} 
B: M \times M & \to \Z/p^n\Z    \\ 
(x_1,y_1), (x_2,y_2) & \mapsto 2(x_1x_2 + \Delta y_1y_2) 
\end{aligned} \]

\noindent on the additive group $G$.
Let $\sigma$ determine a primitive $p^n$th root of unity $\zeta= e^{2\sigma \pi i/p^n}$. The following theorem \cite{nobs1976irreduziblenp1} (pp. 467) allows us to describe a representation in terms of an action on the Bruhat basis. 
\begin{theorem}[Representations] Let $R_k(\Delta, \sigma)$ be the representation on $\C M$ defined by




\begin{equation}\label{eq:a} \left [ \left ( \begin{array}{cc}
a & 0 \\
0 & a^{-1} \\
\end{array} \right ) \right ] \:  \vec{e}_{m}  = \left ( \frac{a}{p} \right )^k {\vec{e}_{m \cdot a^{-1}}} \end{equation}

\begin{equation}\label{eq:b} \left [ \left(
\begin{array}{cc}
1 & b \\
0 & 1 \\
\end{array} \right) \right ] \:   \vec{e}_{m} =  \zeta ^{b  \: Q(m)}  \vec{e}_{m} \end{equation}

\begin{equation}\label{eq:w} \left[ \left(
\begin{array}{cc}
0 & -1 \\
1 & 0 \\
\end{array} \right) \right ] \: \vec{e}_{m} =  c \displaystyle \sum_{m_1 \in M} \zeta^{B(m,m_1)}  \vec{e}_{m_1}  \end{equation}
where $\left ( \frac{a}{p} \right )$ denotes the Legendre symbol and $c$ is a constant 
given by

\[ c = p^{-n +(k/2)} \left ( \frac{\Delta'}{p} \right )^{n -k} \left ( \frac{\sigma}{p} \right )^{k} \cdot e \]
where 

\[e = 
\begin{cases} 1 & k \text{ odd} \hspace{1cm} \left ( \frac{-1}{p} \right ) = 1 \\
-i & k \text{ odd} \hspace{1cm} \left ( \frac{-1}{p} \right ) = -1 \\
-1^n & k  \text{ even} \end{cases}  \]

\end{theorem}

As in the construction described in Section 2, it should be noted that there are three parameters which control the representation $R_k(\Delta, \sigma)$. The integer $k$ controls the dimension of the representations, whereas $\Delta$ and $\sigma$ determine the structure of the representation.  The representations $R_k(\Delta, \sigma)$ are \textbf{not irreducible}, however all irreducible representations of $\slpn$ can be obtained by decomposing the representations $R_k(\Delta, \sigma)$ for various choices of the parameters $k, \Delta, \sigma$. We obtain these subrepresentations by inducing on a character derived from $Q$.

 We endow $G$ with a commutative multiplication structure related to $Q$ $*$ given by 

\[ (x_1,y_1) * (x_2,y_2)  = ( \: x_1x_2 - y_1y_2 \Delta,  \:  x_1y_2 + x_2y_1 \: ) \]

\vspace{.2cm}

While there is an identity element $(1,0) \in G$, this operation is not generally invertible. To understand this definition, let $K = \Q(\sqrt{-\Delta})$ and $\mathcal{O}$ be the order $\mathcal{O} = \Z(\sqrt{-\Delta})$. By identifying $a + b\sqrt{-\Delta}$ with $(a,b)$ then the operation $*$ can be viewed as coming from the standard multiplication of $\mathcal{O}$. Further, the quadratic form $Q$ is just the Norm$_{K\: |\: \Q}$ and the associated bilinear form $B$ is just Tr$_{K\: |\: \Q}(xy^\sigma)$ where $y^\sigma$ is just the non-trivial automorphism in $\text{Gal}( K\: |\: \Q)$ applied to $y$.   Now define

\[ C = \{ \: (x,y) \in G \quad \text{such that }\quad Q(x,y) = 1 \: \} \]

\vspace{.2cm}
\noindent and in fact $(C,*)$ is actually an abelian group. In the literature, $C$ is often defined as the group Aut$(Q,M)$ of automorphisms of $M$ that preserve $Q$. For $c \in C$, then left multiplication by $c$ induces an automorphism of $M$ with $Q(c*m) = Q(m)$ which establishes this correspondence. A classification of when $C$ is cyclic can be found in \cite{nobs1976irreduziblenp1} (pp. 494-495). Further, as in Section 2, the characters on $C$ actually determine subrepresentations of $R_k(\Delta, \sigma)$.  Hence, let $\chi$ be a character on $C$.

\noindent Let

\[ \operatorname{Ind}(V_\chi) :=  \{  \: \vec{v}  \in \C M  \: \: : \: \:  \alpha_{c\: * \: m} = \chi(c) \alpha_{m} \quad \quad c \in C , m \in M \: \}. \]

\vspace{.2cm}
The subspace Ind$(V_\chi)$ is an invariant subspace of $R_k(\Delta, \sigma)$ and we denote this subrepresentation as $R_k(\Delta, \sigma, \chi)$. To work this out more explicitly, let $\sim$ be an equivalence relation on $M$ given by

\[ m_1 \sim m_2 \hspace{1cm} \text{if} \hspace{1cm} m_1 = c*m_2  \hspace{2cm} c \in C. \] 

\vspace{.2cm}
Let $m_i$ be a set of representatives for the equivalence classes $[m_i]$ of this partition. For each equivalence classes $[m_i]$, define

\[ \vec{f_i} = \begin{cases}

0 & \exists c_1 \neq c_2 \in C \text{ with } m_i *c_1=m_i*c_2 \text{ and } \chi(c_1) \neq  \chi(c_2) \\
 \sum_{c \in C} \chi(c)\vec{e}_{c\,* \, m_i}  & \text{Otherwise}  \end{cases}. \]

The vectors $\vec{f_i}$ form a basis for Ind$(V_\chi)$ and the matrices for $R_k(\Delta, \sigma, \chi)$ can be computed by restricting the action of $R_k(\Delta, \sigma)$ to this basis. These are the irreducibles we construct in the main paper. 

\section{Indexing for Experiment Figures}

These tables list  the parameters for $R_k(\chi,\Delta,\sigma)$ for the non--quotient representations of $SL_2(\Z/25\Z)$ in Figures 1-3 and for $SL_2(\Z/125\Z)$ in Figure 4. The index column is the ordering used along the y-axis in these figures. 

\begin{table}[!h]
\begin{tabular}{|c|c|c|c|c||c|c|c|c|c|}
\hline
Index&$k$&$\chi$&$\Delta$&$\sigma$&Index&$k$&$\chi$&$\Delta$&$\sigma$\\
\hline
\hline
1&0&1&1&1    & 21 &1&1&1&1\\
\hline
2&0&2&1&1    & 22 &1&2&1&1\\
\hline
3&0&3&1&1    & 23 &1&3&1&1\\
\hline
4&0&4&1&1    & 24 &1&4&1&1\\
\hline
5&0&6&1&1 & 25 &1&1&2&1\\ 
\hline
6&0&7&1&1 & 26 &1&2&2&1\\ 
\hline
7&0&8&1&1 & 27 &1&3&2&1\\ 
\hline
8&0&9&1&1 & 28 &1&4&2&1\\ 
\hline
9&0&11&1&1 & 29 &1&1&1&2\\ 
\hline
10&0&12&1&1 & 30 &1&2&1&2\\ 
\hline 
11&0&13&1&1 & 31 &1&3&1&2\\ 
\hline 
12&0&14&1&1 & 32 &1&4&1&2\\ 
\hline 
13&0&1&2&1  & 33 &1&1&2&2\\ 
\hline 
14&0&2&2&1 & 34 &1&2&2&2\\ 
\hline 
15&0&3&2&1 & 35 &1&3&2&2\\ 
\hline
16&0&4&2&1 & 36 &1&4&2&2\\ 
\hline
17&0&6&2&1& 37&2&1&1&1\\ 
\hline
18&0&7&2&1& 38&2&-1&1&1\\ 
\hline
19&0&8&2&1& 39&2&1&1&2\\ 
\hline
20&0&9&2&1 & 40&2&-1&1&2\\ 
\hline
\end{tabular}
\caption{Indices for y--axis in Figures 1--3.}
\end{table}

\begin{table}[!h]
\begin{tabular}{|c|c|c|c|c|}
\hline
Index&$k$&$\chi$&$\Delta$&$\sigma$\\
\hline
\hline
1&1&1&1&1\\
\hline
2&1&2&1&1\\
\hline

3&2&1&1&1\\
\hline

4&2&2&1&1\\
\hline

\end{tabular}
\caption{Indices for y--axis in Figures 1--3.}
\end{table}

\end{document}